%AMSTeX File

% masl15.tex --- 17-01-96 12:20

\documentstyle{amsppt}
\TagsOnRight
%\NoBlackBoxes
\hoffset=.1in

\define\du{\partial}
\define\duo{\partial_0}
\define\dui{\partial_i}
\define\dum{\partial_\mu}
\define\baba{e^{i \frac{\duo}{\kappa}}}
\define\kot{\dfrac{\partial f}{\partial x^i}}
\define\kotek{\dfrac{\partial f}{\partial x^\mu}}

\define\lx{\operatorname{L}}

\define\id{\operatorname{id}}

\define\lrar{\Longleftrightarrow }

\define\duch{\chi_{{\phantom{}}_\mu}}   %?=indeks
\define\haduchmu{\widehat{\chi}_{{\phantom{}}_\mu}}
\define\hachi{\widehat{\chi}}

\define\a{\alpha}
\define\be{\beta}

\define\ve{\varepsilon}

\define\de{\delta}

\define\vd{\varDelta}
\define\vl{\varLambda}

\define\cam{\Cal M}

\define\rti{\widetilde \rho}

\define\hac{\widehat C}
\define\ham{\widehat M}
\define\han{\widehat N}
\define\hap{\widehat P}
\define\hax{\widehat X}

\define\pd#1#2{\dfrac{\partial#1}{\partial#2}}

\define\vc#1{(#1_1,\ldots,#1_n)}
\define\vct#1{[#1_1,\ldots,#1_n]}
\define\vect#1{\{#1_1,\ldots,#1_n\}}

\define\ft{\infty}

\define\xx{\otimes}

\define\jak{\dfrac{i}{\kappa}}

\define\maks{{P_0}\slash{\kappa}}
\define\vep{\vec{P\,}^2}

%\define\chi_{{\phantom{}}_{ \chi_{{\phantom{}}_} }

\define\rol{\rho_{\lx}}

\define\tirol{\rti_{\lx}}

\define\lama{\vl^\mu{}_\nu}

\define\kap{{\Cal P}_\kappa}
\define\kapti{\widetilde{\Cal P}_\kappa}
\define\mina{\cam_\kappa}

\define\1{$\kappa$-Poincar\'e group}
\define\2{covariant}
\define\3{calculus}
\define\4{differential}
\define\5{dimensional}
\define\6{quantum}
\define\7{quantum group}
\define\8{condition}
\define\9{generalization}
\define\0{algebra}

\define\defx{deformation}
\define\pro{properties}
\define\wor{Wo\-ro\-no\-wicz}

\define\fol{following}

\define\kapo{$\kappa$-Poincar\'e}
\define\mink{$\kappa$-Min\-ko\-wski}
\define\linv{left-covariant}

\define\five{five-dimensional}
\define\fif{fifteen-dimen\-sion\-al}
\define\ina{infinitesimal action}

\define\ajka{a_{(1)}}
\define\adka{a_{(2)}}
\define\xaj{x_{(1)}}
\define\xad{x_{(2)}}
\define\ajd{a_{(1)(2)}}

\topmatter
\title A note on geometry of $\kappa$-Minkowski space
\endtitle
\rightheadtext{Differential calculi}
\author Stefan Giller$^*$, Cezary Gonera$^*$, Piotr Kosi\'nski$^*$, Pawe\l 
\/  Ma\'slanka$^*$\\
{\it Department of Theoretical Physics}\\
{\it University of \L \'od\'z}\\
{\it ul. Pomorska 149/153, 90--236 \L \'od\'z, Poland}
\endauthor
\leftheadtext{S. Giller, C. Gonera, P. Kosi\'nski, P. Ma\'slanka}
\thanks
*\ \ \ Supported by \L \'od\'z University grant N$^o$ 458
\endthanks

\abstract The infinitesimal action of $\kappa$-Poincar\'e group on 
$\kappa$-Minkowski space is computed both for generators of 
$\kappa$-Poincar\'e algebra and those of \wor{} generalized Lie \0. The 
notion of invariant operators is introduced and generalized Klein-Gordon 
equation is written out.
\endabstract

\endtopmatter

\document

\head I. Introduction
\endhead

In this short note we consider some simple \pro{} of \4 operators on  
\mink{} space $\mina$ --- a noncommutative \defx{} of Min\-kowski space-time 
which depends on dimensionful parameter $\kappa$ ([1]). We  calculate the 
infinitesimal action  of $\kappa$-Poincar\'e group  $\kap$ ([1]) on 
$\mina$  both for the generators of $\kappa$-Poincar\'e algebra $\kapti$ 
([2]) (this is done using the duality $\kapti \lrar \kap$ described in [3]) 
and for the elements of \wor{} generalized Lie \0 ([4]) of 
$\kappa$-Poincar\'e group  ([5]). The result supports the relation between 
both \0s found in [5]. We introduce also the notion of invariant \4 
operators on $\mina$ and write out the generalized Klein-Gordon equation.

Let us conclude this section by introducing the notions of \1 $\kap$ and \0 
$\kapti$.  \ $\kap$ is defined by the \fol{} relations ([1])
$$
\aligned
&[x^\mu,x^\nu]  = \dfrac{i}{\kappa}( \de_0^\mu x^\nu -  \de_0^\nu 
 x^\mu),\\
&[\lama,  \vl^\a{}_\be]  = 0,\\
& [\lama,  x^\rho]  = -  \dfrac{i}{\kappa} ((\vl^\mu{}_0 - \de_0^\mu) 
\vl^\rho{}_\nu + (\vl^0{}_\nu - \de_\nu^0)  g^{\mu\rho}),\\
& \vd(\lama)  = \vl^\mu{}_\a \otimes \vl^\a{}_\nu,\\
& \vd(x^\mu)  = \vl^\mu{}_\a \otimes x^\a + x^\mu  \otimes I,\\
& S(\vl^\mu{}_\nu) =  \vl_\nu{}^\mu,\\
& S(x^\mu) = - \vl_\nu{}^\mu x^\nu,\\ 
& \ve(\vl^\mu{}_\nu)  =  \de^\mu_\nu,\\
& \ve(x^\mu) = 0.
\endaligned
\tag{1}
$$
The dual structure, the \kapo{} \0 $\kapti$, is, in turn, defined as follows 
([6])
$$
\aligned
& [P_\mu,P_\nu] = 0,\\
& [M_i,M_j] = i \ve_{ijk} M_k,\\
& [M_i,N_j] = i \ve_{ijk} N_k,\\
& [N_i,N_j] = - i \ve_{ijk} M_k,\\
& [M_i,P_0] = 0,\\
& [M_i,P_j] = i\ve_{ijk}P_k,\\
& [N_i,P_0] = iP_i,\\
& [N_i,P_j] = i \de_{ij}  \Big( \dfrac{\kappa}{2}(1 - e^{-2\maks}) + 
\dfrac{1}{2\kappa} \vep \Big) - \dfrac{i}{\kappa} P_iP_j,\\
& \vd (M_i) =  M_i  \otimes I + I \otimes M_i,\\
& \vd (N_i)  = N_i  \otimes e^{-\maks} + I \otimes N_i - 
\dfrac{1}{\kappa} \ve_{ijk} M_j \otimes  P_k,\\
& \vd (P_0)  = P_0 \otimes I + I \otimes P_0,\\
& \vd (P_i) = P_i \otimes e^{-\maks} + I \otimes P_i,\\
&S(M_i) = -M_i,\\
&S(N_i) = - N_i + \dfrac{3i}{2\kappa} P_i ,\\
&S(P_\mu) = - P_\mu,\\ 
& \ve(P_{\mu},M_i,N_i) = 0.
\endaligned
\tag{2}
$$

Structures (1), (2) are dual to each other, the duality being fully 
described in [3].

The analysis given below was suggested to two of the  authors (P. Kosi\'nski 
and P. Ma\'slanka) by J. Lukierski.

\head II. $\kappa$-Min\-ko\-wski space
\endhead

The $\kappa$-Min\-ko\-wski space  $\mina$ ([1]) is a universal $\ast$-\0 
with unity generated by four selfadjoint elements $x^\mu$ subject to the 
\fol{} \8s
$$
[x^\mu, x^\nu] =  \dfrac{i}{\kappa} (\de^\mu_0 x^\nu - \de^\nu_0 x^\mu ) .
\tag"{(3a)}"
$$
Equipped with the standard coproduct
$$
\vd x^\mu =  x^\mu \otimes  I + I \otimes  x^\mu ,
\tag"{(3b)}"
$$
antipode $S(x^\mu) = - x^\mu$ and counit $\ve(x^\mu) = 0$ it becomes a \7.

On $\mina$ one can construct a bi\2 \five{} \3 which is defined by the 
\fol{} relations ([5])
$$
\aligned
& \tau^\mu \equiv d x^\mu, \qquad \qquad \tau \equiv d\Big( x^2 + 
\dfrac{3i}{\kappa}x^0 \Big) - 2 x_\mu  d x^\mu,\\
& [\tau^\mu, x^\nu] = \jak g^{0\mu} \tau^{\nu} - \jak g^{\mu\nu} \tau^0 + 
\dfrac{1}{4}  g^{\mu\nu} \tau,\\
& [\tau,x^\mu]  = - \dfrac{4}{\kappa^2}  \tau^{\mu},\\
& \tau^\mu \wedge  \tau^\nu =  -\tau^\nu \wedge \tau^\mu     ,\\
& \tau \wedge \tau^\mu   = -  \tau^\mu \wedge  \tau,\\
& (\tau^\mu)^*  =  \tau^\mu, \qquad\qquad \tau^* = -\tau,\\
& d\tau^\mu = 0,\\
& d\tau = - 2\tau_\mu  \wedge \tau^\mu .
\endaligned
\tag{4}
$$
The \mink{}  space carries a \linv{} action of \1   $\kap$ ([1]), $\rol : 
\mina \to \kap \otimes \mina$,  given by
$$
\rol (x^\mu) = \lama  \xx x^\nu + a^\mu \xx I.
\tag{5}
$$

The \3 defined by (4) is \2 under the action of $\kap$ which reads
$$
\aligned
& \tirol (\tau^\mu) = \lama \xx \tau^\nu, \\
& \tirol (\tau) = I \xx \tau.
\endaligned
\tag{6}
$$

\head III. Derivatives, infinitesimal actions and invariant operators 
\endhead

The product of generators $x^\mu$ will be called normally ordered if all 
$x^0$ factors stand leftmost. This definition can be used to ascribe a 
unique element $:f(x) :$ of $\mina$ to any polynomial function of four 
variables $f$. Formally, it can be extended to any analytic function $f$.

Let us now one define the (left) partial derivatives: for any $f \in \mina$ 
we write
$$
df = \partial_\mu f \tau^\mu + \partial f \tau.
\tag{7}
$$

It is a matter of some boring calculations (using the commutation rules 
(3a)) to find the \fol{} formula
$$
\aligned
d:f: & =  : \Big( \kappa \sin \big(\dfrac{\duo}{\kappa} \big) + \dfrac{i}{2\kappa} 
\baba \vd\Big)f : \tau^0 + : \baba \kot : \tau^i\\
& + : \Big( \dfrac{\kappa^2}{4} \Big(1 - \cos \big(\dfrac{\duo}{\kappa} 
\big) \Big) - \dfrac{1}{8} \baba  \vd\Big)f : \tau
\endaligned
\tag{8}
$$
or
$$
\aligned
& \duo :f:  =  : \Big(\kappa \sin \big(\dfrac{\duo}{\kappa} \big) + 
\dfrac{i}{2\kappa} \baba \vd\Big)f : \\
& \dui :f:  =  :  \baba \kot: \\
& \du : f : = : \Big( \dfrac{\kappa^2}{4} \Big(1 - \cos 
\big(\dfrac{\duo}{\kappa}  \big) \Big) - \dfrac{1}{8} \baba  \vd\Big)f :
\endaligned
\tag{9}
$$

Let us now define the infinitesimal action of $\kap$ on $\mina$. Let $X$ be 
any element of the Hopf \0 dual to $\kap$ --- the \kapo{} \0 $\kapti$ (cf. 
[3] for the proof of duality). The corresponding \ina{}
$$
\hax : \mina \to \mina
$$ 
is defined as follows: for any $f \in \mina$,
$$
\hax f = (X \xx \id) \circ \rol (f).
\tag{10}
$$
Using the standard duality rules ([3]), we conclude that
$$
\aligned
& \hap_\mu x^\a = i \de_\mu^\a,\\
& \hap_\mu : x^\a  x^\be : = i  \de_\mu^\be x^\a + i \de_\mu^\a x^\be
\endaligned
\tag{11}
$$
etc. One can show that, in general,
$$
 \hap_\mu : f : = : i \kotek :
\tag{12}
$$

Also, using the fact that $\tirol$ is a left action of $\kap$ on $\mina$ 
together with the duality $\kap \to \kapti$, we conclude that
$$
F( \hap_\mu) : f : = : F \Big( i \kotek  \Big) f :
\tag{13}
$$

Formulae (11)--(13) have the \fol{} interpretation. In [5] the \fif{} bi\2 
\3 on $\kap$ has been constructed using the methods developed by \wor{} 
([4]).  The resulting generalized Lie \0 is also \fif{}, the additional 
generators being the generalized mass square operator and  the components of 
generalized Pauli-Lubanski fourvector. All generators of this Lie \0 can be 
expresses in terms of the generators $P_\mu$, $M_{\a\be}$ of $\kapti $ 
([5]). In particular, the translation generators $\duch$ as well as the mass 
squared  operator $\chi$ are expressible in terms of $P_\mu$ only. The 
relevant expressions are given by formulae (20) of [5]. Comparing them with 
(9), (13) above, we conclude that
$$
\aligned
& \haduchmu \equiv \dum,\\
& \hachi \equiv \du.
\endaligned
\tag{14}
$$

These relations, obtained here by explicit computations, follow also from 
(7) if one takes into account that $\mina$ is a \6 subgroup of $\kap$.

It is also not difficult to obtain the action of Lorentz generators. 
Combining (1) and (3a) with the  duality  $\kap \to \kapti$ described in 
detail in [5], we conclude first that the action of $M_i$ and $N_i$ 
coincides with the proposal of Majid and Ruegg ([6]); the actual computation 
is then easy and gives
$$
\aligned
& \ham_i : f(x^\mu) : = : -i \ve_{ijl} x^j \pd{f(x^\mu)}{x^l}:\\
& \han_i : f(x^\mu) : = : \Big(i x^0 \pd{}{x^i} + x^i \Big(\dfrac{\kappa}{2} 
\Big(1 - e^{\frac{-2i}{\kappa}\frac{\du}{\du x^0}}\Big) - \dfrac{1}{2\kappa} 
\vd\Big)\\
& \phantom{\han_i : f(x^\mu) : = :\qquad}  + \dfrac{1}{\kappa} x^k 
\dfrac{\du^2}{\du x^k \du x^i} \Big) f(x^\mu):
\endaligned
\tag{15}
$$

Let us now pass to the notion of invariant operator; $\hac$ is an invariant 
operator on $\mina$ if
$$
\rol \circ \hac = (\id \xx \hac ) \circ \rol.
\tag{16}
$$

We shall show that  if $C$ is a central element of $\kapti$, then
$$
\hac f = (C \xx \id) \circ \rol(f)
\tag{17}
$$
is an invariant operator. To prove this let us take any $Y \in \kapti$, then
$$
YC = CY
\tag{18}
$$
or, in other words, for any $a \in \kap$,
$$
Y(\ajka) C(\adka) = C(\ajka)Y(\adka)
\tag{19}
$$
where $\vd a = \ajka \xx \adka$. Let us fix $a$ and write (19) as
$$
Y(\ajka C(\adka)) = Y(\adka C(\ajka)).
\tag{20}
$$
As (20) holds for any $Y \in \kapti$ we conclude that for any $a \in \kap$
$$
C(\adka) \ajka = C(\ajka) \adka.
\tag{21}
$$

Now let
$$
\aligned
& \rol (x) = \ajka \xx \xaj,\\
& \rol(\xaj) = a_{(1)(2)} \xx \xad,\\
& \vd \ajka = a_{(1)}^{(1)} \xx a_{(1)}^{(2)}.
\endaligned
\tag{22}
$$

The identity
$$
 (\id \xx \rol ) \circ \rol = (\vd \xx \id ) \circ \rol
\tag{23}
$$
implies
$$
\ajka \xx \ajd \xx \xad = a_{(1)}^{(1)} \xx a_{(1)}^{(2)} \xx \xaj.
\tag{24}
$$

Applying to both sides $\id \xx C \xx \id$ and $C \xx \id \xx \id$, we get
$$
\aligned
& C(\ajd) \ajka \xx \xad = C(a_{(1)}^{(2)}) a_{(1)}^{(1)} \xx \xaj,\\
& C(\ajka) \ajd \xx \xad = C(a_{(1)}^{(1)}) a_{(1)}^{(2)} \xx \xaj.
\endaligned
\tag{25}
$$

%a_{(1)}^{(1)} \xx a_{(1)}^{(2)}

It follows from (21) applied to $\ajka$ that the right-hand sides of (25) 
are equal. So,
$$
C(\ajd) \ajka \xx \xad = C( \ajka) \ajd  \xx \xad
\tag{26}
$$
i.e.
$$
 (\id \xx \hac) \circ \rol(x)  =  \rol \circ \hac(x).
\tag{27}
$$

Using the above result we can easily construct the deformed Klein-Gordon 
equation. Namely, we take as a central element the counterpart of mass 
squared Casimir operator $\chi$ ([5]). Due to (14) the generalized  
Klein-Gordon equation reads
$$
\Big( \du + \dfrac{m^2}{8}\Big) f = 0;
\tag{28}
$$
the coefficient $\frac{1}{8}$ is dictated by the correspondence with 
standard  Klein-Gordon equation in the limit $\kappa \to \ft$. Let us note 
that (28) can be written, due to (9), in the form
$$
\Big[ \du_0^2 - \du_i^2 + m^2 \Big( 1 + \dfrac{m^2}{4\kappa^2}\Big) \Big] f 
= 0; 
\tag{29}
$$
here $\du_0$, $\du_i$ are the operators given by (9). It seems therefore 
that the \wor{} operators $\duch$ are better candidates for translation 
generators than $P_\mu$'s. Note that the operators $\duch$ already appeared  
in [7], [8].


\Refs
\ref \key 1 \by S. Zakrzewski  \jour  J. Phys. \vol A\,27 \yr 1994 \pages 
2075 \endref
\ref \key 2 \by J. Lukierski, A. Nowicki, H. Ruegg \jour  Phys. Lett. \vol 
B\,293 \pages  344 \yr 1992  \endref
\ref\key 3 \by  P. Kosi\'nski,  P. Ma\'slanka \paper The duality between 
$\kappa$-Poincar\'e algebra and $\kappa$-Poincar\'e group \jour hep-th 
\endref 
\ref \key 4 \by S.L. Woronowicz \jour Comm. Math. Phys. \vol 122 \pages  125 
\yr  1989 \endref
\ref\key 5 \by   P. Kosi\'nski,  P. Ma\'slanka, J. Sobczyk \paper in: Proc. 
of IV Coll. on Quantum Groups and Integrable Sysytems, Prague 1995, to be 
published  \endref
\ref\key 6 \by S. Majid,  H. Ruegg \jour  Phys. Lett. \vol B\,334 \pages  
348 \yr 1994  \endref
\ref\key 7 \by  H. Ruegg, V. Tolstoy \jour Lett. Math. Phys. \vol 32 \pages 
85 \yr  1994 \endref
\ref\key 8 \by J. Lukierski,  H. Ruegg, V. Tolstoy \jour in: Quantum Groups. 
Formalism and Applications.  Proc. of XXX Karpacz Winter School of 
Theoretical  Physics, ed. J. Lukierski, Z. Popowicz, J. Sobczyk,  PWN 1995 
\endref

\endRefs
\enddocument